\documentclass[11pt,a4paper]{amsart}
\pdfoutput=1

\usepackage[utf8]{inputenc}
\newcommand{\boldm}[1] {\mathversion{bold}#1\mathversion{normal}}
\usepackage{hyperref,graphicx,amssymb}
\newtheorem*{theorem*}{Theorem}

\begin{document}
\title[Bounding the transitivity degree of Galois groups]{On Elkies' method for bounding the transitivity degree of Galois groups}

\author{Dominik Barth}
\author{Andreas Wenz}

\address{Institute of Mathematics\\ University of Würzburg \\ Emil-Fischer-Straße 30 \\ 97074 Würzburg, Germany}
\email{dominik.barth@mathematik.uni-wuerzburg.de}
\email{andreas.wenz@mathematik.uni-wuerzburg.de}

\begin{abstract}
In 2013 Elkies described a method for bounding the transitivity degree of Galois groups. Our goal is to give additional applications of this technique, in particular verifying that the monodromy group of the degree-$276$ cover defined over a degree-$12$ number field computed by Monien is isomorphic to the sporadic Conway group $\text{Co}_3$.
\end{abstract}

\maketitle

\section{Introduction}

 With the recent development of computing explicit polynomials of large degree with prescribed Galois groups the corresponding verification process poses new computational challenges. 
 
 While standard methods, for example Dedekind criterion for obtaining lower bounds for Galois groups, are still rather viable for polynomials of large degree, other techniques, such as the calculation of resolvents, are not expected to be feasible. In particular, the question whether a $2$-transitive Galois group is a full symmetric or alternating group is generally difficult to answer.

 In 2013 Elkies \cite{Elkies} described a method for bounding the transitivity degree of Galois groups of function field extensions
by collecting the factorization patterns of many specialized polynomials and comparing them to an effective version of Chebotarev's density theorem which arises from the Hasse-Weil bound.
Elkies used this particular technique to verify that the Galois group of his computed polynomial is the sporadic $4$-transitive Mathieu group $M_{23}$.

 Our goal is to give additional applications of this technique, in particular rigorously verifying that the monodromy group of the degree-$276$ cover defined over a degree-$12$ number field computed by Monien \cite{Monien} is isomorphic to the sporadic $2$-transitive Conway group $\text{Co}_3$.
 
 This paper is structured as follows: section \ref{preliminaries} introduces the objects of interest, section \ref{methodbyelkies} depicts Elkies' technique for bounding the transitivity degree of Galois groups. In the final section \ref{newapplications} we present several new applications of Elkies' method dealing with the Galois groups $\text{Co}_3$, $\text{PSp}_6(2)$ and $\text{PSL}_6(2)$.

\section{Preliminaries}\label{preliminaries}
For a fixed number field $K$ let $p$ and $q$ be coprime polynomials in $K[X]$. The arithmetic monodromy group of the degree-$n$ cover $\frac{p}{q}$ is defined as 
\[
A:= \text{Gal}(N \mid K(t))
\]
where $N$ denotes the splitting field of $p(X)-tq(X)$ over $K(t)$. Furthermore, the geometric monodromy group of $\frac{p}{q}$ is defined as
\[
G := \text{Gal}(N \mid (\bar{K} \cap N)(t)).
\] 
Since $p(X)-tq(X)$ is absolutely irreducible the natural (faithful) action of both $G$ and $A$ on the $n$ roots of $p(X)-tq(X)$ in $N$ is transitive. Furthermore, it is well known that $G$ is normal in $A$.

In order to study $A$ and $G$, we will reduce the above polynomials modulo a suitable prime: 
The ring of integers of $K$ will be denoted by $\mathcal{O}_K$.
For a fixed prime ideal $\mathfrak{p}$ in $\mathcal{O}_K$ we write $p_\mathfrak{p}$ and $q_\mathfrak{p}$ for the reduction of $p$ and $q$ modulo $\mathfrak{p}$. In the same fashion as before we define 
\begin{align*}
 A_\mathfrak{p} := \text{Gal} (N_\mathfrak{p} \mid (\mathcal{O}_K/\mathfrak{p})(t))
\quad \text{and} \quad 
 G_\mathfrak{p} := \text{Gal} (N_\mathfrak{p}\mid (\overline{\mathcal{O}_K/\mathfrak{p}}\cap N_\mathfrak{p})(t))
\end{align*}
 where  $N_\mathfrak{p}$ denotes the splitting field of $p_\mathfrak{p}(X)- t q_\mathfrak{p} (X)$ over $(\mathcal{O}_K/\mathfrak{p})(t)$. Again, $G_\mathfrak{p}$ is a normal subgroup of $A_\mathfrak{p}$.

Thanks to a theorem of Beckmann \cite[Proposition 10.9]{MM}, among other considerations, if $\mathfrak{p}$ is chosen to be lying over a sufficiently large (rational) prime we may assume the following: 
\begin{enumerate}
\item[(i)] The ramification locus of $p(X)-tq(X)$ with respect to $t$ is $\mathfrak{p}$-stable.
\item[(ii)] The inseparability behaviour of both $p(X)-tq(X)\in K(t)[X]$ and $p_\mathfrak{p}(X)-tq_\mathfrak{p}(X)\in (\mathcal{O}_K/\mathfrak{p})(t)[X]$ specialized at ramified places with respect to $t$ coincides.
\item[(iii)] $G \cong G_\mathfrak{p}$.
\end{enumerate}

\section{A method by Elkies} \label{methodbyelkies}

The following technique described by Elkies (see \cite{Elkies}) bounds the transitivity degree of $G$: 

Assume, $G$ and therefore $G_\mathfrak{p}$ is $k$-transitive and $A_\mathfrak{p} = G_\mathfrak{p}$. Let $C_0$ and $C_1$ be the projective $t$- and $x$-lines over the finite field $\mathbb{F}_\lambda \cong \mathcal{O}_K/\mathfrak{p}$. By introducing the relation $p_\mathfrak{p}(x)-tq_\mathfrak{p}(x)=0$ we obtain a cover $C_1/C_0$ ramified over exactly $m$ points with ramification structure $(s_1,\dots,s_m) \in (S_n)^m$. Its Galois closure will be denoted by $\tilde{C}$.

Let $(G_\mathfrak{p})_k$ be the stabilizer of a $k$-element set in $G_\mathfrak{p}$ and $C_k := \tilde{C}/(G_\mathfrak{p})_k$. 
The corresponding cover $C_k/C_0$ is of degree $\binom{n}{k}$ with ramification structure $(\sigma_1,\dots,\sigma_m)$ induced by the natural action of $(s_1,\dots,s_m)$ on $k$-element subsets. As $G_\mathfrak{p}$ acts faithfully on $n$ elements, it can be shown easily that the action on $k$-element subsets is also faithful if $k\not\in \lbrace 0,n \rbrace$. In particular, $\text{ord}(\sigma_i) = \text{ord}(s_i)$ for $i=1, \dots ,m$.
Additionally note that $C_k$ is an irreducible curve with full constant field $\mathbb{F}_\lambda$ due to $A_\mathfrak{p}=G_\mathfrak{p}$.

The number of $\mathbb{F}_\lambda$-rational points on $C_k$, denoted by $\#C_k(\mathbb{F}_\lambda)$, has to obey the Hasse-Weil bound $\left| \#C_k(\mathbb{F}_\lambda) - (\lambda + 1) \right| \leq 2g(C_k)\sqrt{\lambda}$, in particular
\begin{equation}
\label{hasse-weil}
\#C_k(\mathbb{F}_\lambda) \leq \lambda + 1 + 2g(C_k)\sqrt{\lambda}.
\end{equation}
Here, $g(C_k)$ denotes the genus of $C_k$.
In order to check if $C_k$ is indeed compatible with the above bound, we need to determine $\#C_k(\mathbb{F}_\lambda)$ and $g(C_k)$.

We will use the following notation: For a permutation $s\in S_n$ let $\pi_k(s)$ be the number of invariant $k$-element subsets of $s$.

\subsection{Counting {\boldm $\mathbb{F}_\lambda$-rational points on $C_k$}}

Fix $t_0\in \mathbb{P}^1(\mathbb{F}_\lambda)$ not contained in the ramification locus $S$ of $p_\mathfrak{p}(X)-tq_\mathfrak{p}(X)$. 
Note that $\mathbb{F}_\lambda$-rational points on $C_k$ lying over $t_0$ correspond to degree-$k$ factors of the specialization $p_\mathfrak{p}(X)-t_0q_\mathfrak{p}(X)\in \mathbb{F}_\lambda[X]$. 

If $\text{Frob}(t_0)$ denotes the Frobenius permutation on the $n$ roots of the specialization $p_\mathfrak{p}(X)-t_0q_\mathfrak{p}(X)$
then the number of $\mathbb{F}_\lambda$-rational points on $C_k$ lying over $t_0$ is given by $\pi_k(\text{Frob}(t_0))$, therefore 
\begin{equation}\label{zaehlen}
\#C_k(\mathbb{F}_\lambda) \geq \sum_{t_0 \in \mathbb{P}^1(\mathbb{F}_\lambda) \setminus S} \pi_k(\text{Frob}(t_0)).
\end{equation}

\subsection{Computing the genus of {\boldm $C_k$}}
Since the degree-$\binom{n}{k}$ cover $C_k/C_0$ has ramification structure $(\sigma_1,\dots,\sigma_m)$ the Riemann-Hurwitz formula yields
\begin{equation}\label{g}
g(C_k) = 1-\binom{n}{k} + \frac{1}{2} \sum_{i=1}^m \text{ind}(\sigma_i)
\end{equation}
where $\text{ind}(\sigma_i):= \binom{n}{k} - \text{number\;of\;cycles\;of\;}\sigma_i$. 

If $(\sigma_1,\dots,\sigma_m)$ cannot be computed explicitly, one can deduce the upper bound 
\begin{equation}\label{ind}
\text{ind}(\sigma_i) \leq \left( \binom{n}{k}-\pi_k(s_i) \right) \left( 1-\frac{1}{\text{ord}(s_i)} \right).
\end{equation}
Note that equality holds if the order of $s_i$ is prime.

\subsection{Picking a sufficiently large prime}

Note that the right hand side of \eqref{zaehlen} behaves differently if $G_\mathfrak{p}$ is not $k$-transitive: Let $d$ be the number of orbits of $G_\mathfrak{p}$ acting on $k$-element subsets, then it is reasonable to expect
\begin{equation}\label{cheb}
\sum_{t_0\in \mathbb{P}^1(\mathbb{F}_\lambda)\setminus S} \pi_k(\text{Frob}(t_0)) \approx d \lambda
\end{equation} 
for large $\lambda$ due to the orbit-counting theorem in combination with Chebotarev's density theorem.
By comparing \eqref{zaehlen} and \eqref{cheb} with the Hasse-Weil bound \eqref{hasse-weil} we obtain $d\lambda \leq \lambda + 2g(C_k)\sqrt{\lambda}$, which leads to
$\lambda \leq \frac{4g(C_k)^2}{(d-1)^2}$ in the case $d>1$. 

This observation is a crucial ingredient in the verification process: If $\mathfrak{p}$ is picked such that its norm $\lambda$ is sufficiently greater than $\frac{4g(C_k)^2}{(d-1)^2}$, we are able to establish a contradiction to the Hasse-Weil bound. This, in particular, allows us to disprove the $k$-transitivity of $G$.

\subsection{Elkies' example: {\boldm $M_{23}$}} \label{M} The previously described technique was key for the proof that the degree-$23$ polynomial $p\in \mathbb{Q}(\alpha)[X]$ given in \cite{Elkies} where $\alpha^4+\alpha^3+9\alpha^2-10\alpha+8=0$ has geometric monodromy group $M_{23}$. 

With respect to $t$ the ramification locus of $p(X)-t$ consists of exactly three points with ramification type $(4^4.2^2.1^3, 2^8.1^7, 23^1)$.
It is standard practice to show that the geometric monodromy group $G$ of $p$ is either the $4$-transitive group $M_{23}$ or $A_{23}$. 

Assume $G \cong A_{23}$, then for the prime ideal $\mathfrak{p}:=(47\,000\,081, \alpha + 25\,037\,440)$ of norm $\lambda := 47\,000\,081$ we have $G\cong G_\mathfrak{p}$ and both $G$ and $G_\mathfrak{p}$ are $5$-transitive. Since the discriminant of $p_\mathfrak{p}(X)-t$ is a square, $A_\mathfrak{p} \cong A_{23} \cong G_\mathfrak{p}$.
This leads us to work with the curve $C_5$: 
By explicitly computing $(\sigma_1,\sigma_2,\sigma_3)$ we find $g(C_5)=3285$ using the Riemann-Hurwitz formula \eqref{g}. In combination with the Hasse-Weil bound \eqref{hasse-weil} this yields $\#C_5(\mathbb{F}_\lambda) \leq 92\,041\,771$. Counting $\mathbb{F}_\lambda$-rational points on $C_5$ according to \eqref{zaehlen} reveals the contradiction $\#C_5(\mathbb{F}_\lambda) \geq 93\,981\,891$ with a total computing time of approximately 12 hours using \texttt{Magma} \cite{Magma}. We obtain $G\cong G_\mathfrak{p}\cong M_{23}$.

\section{New Applications} \label{newapplications}

\subsection{The sporadic Conway group {\boldm $\text{Co}_3$}} \label{Co}
In this section we will refer to the polynomials $p:=-k_3 \tilde{p}_3$ and $q:=k_2 \tilde{p}_2$ presented in \cite[Proposition 1]{Monien} of degree $276$ over a degree-$12$ number field $K:=\mathbb{Q}(\alpha)$ where
$
\alpha^{12}-
2\alpha^{11} + 
9\alpha^{10} -
20\alpha^9 + 
38\alpha^8 - 
73\alpha^7 + 
101\alpha^6 - 
86\alpha^5 + 
55\alpha^4 - 
46\alpha^3 + 
42\alpha^2 - 
24\alpha + 
6=0
$.

\begin{theorem*}
The polynomial $p(X)-tq(X)\in K(t)[X]$ defines a regular Galois extension of $K(t)$ with Galois group isomorphic to the sporadic $2$-transitive Conway group $\mathrm{Co}_3$. With respect to $t$ the ramification locus is given by $(0,1,\infty)$ with corresponding ramification type $(3^{92},7^{39}.1^3,2^{132}.1^{12})$.
\end{theorem*}
\begin{proof}
An easy computation shows that $p(X)-tq(X)$ is ramified over 0, 1 and $\infty$ with the given ramification type. The ramification locus cannot be any larger, otherwise this would contradict the Riemann-Hurwitz formula.

Pick the prime ideal $\mathfrak{p}:=(7 \cdot 10^{9} +1,\alpha + 2\,738\,443\,742)$ in $\mathcal{O}_{K}$ of norm $\lambda := 7 \cdot 10^{9} +1$. Note that $\mathcal{O}_K /\mathfrak{p} \cong\mathbb{F}_\lambda$. Because $\frac{p_\mathfrak{p}(X)q_\mathfrak{p}(t)-p_\mathfrak{p}(t)q_\mathfrak{p}(X)}{X-t} \in \mathbb{F}_\lambda(t)[X]$ is irreducible, $A_\mathfrak{p}$ must be $2$-transitive. Additionally, the discriminant of $p_\mathfrak{p}(X)-tq_\mathfrak{p}(X)\in \mathbb{F}_\lambda(t)[X]$ is a square. Combining both results, we find $A_\mathfrak{p} \in \{\text{Co}_3, A_{276}\}$ by the classification of finite $2$-transitive groups. In both cases we have $G_\mathfrak{p} = A_\mathfrak{p}$ because $G_\mathfrak{p}$ is normal in $A_\mathfrak{p}$. 

Under the assumption that $G_\mathfrak{p}$ is $3$-transitive we study the curve $C_3$:
Combining $\eqref{ind}$ and $\eqref{g}$ yields $g(C_3) = 40\,782$.
Now, \eqref{zaehlen} gives us
$
\#C_3(\mathbb{F}_\lambda) \geq 13\,999\,925\,705
$
whereas $\#C_3(\mathbb{F}_\lambda)\leq 13\,824\,133\,842$ by the Hasse-Weil bound \eqref{hasse-weil}. This is a contradiction, thus
$G_\mathfrak{p}$ cannot be $3$-transitive and we remain with $G_\mathfrak{p} = \text{Co}_3$. 

Since $\mathfrak{p}$ is a prime of good reduction for $p(X)-tq(X)\in K(t)[X]$, a theorem of Beckmann, see \cite[Proposition 10.9]{MM}, implies $G\cong G_\mathfrak{p} = \text{Co}_3$. Due to the fact that $G$ is normal in $A$ and $N_{S_{276}}(\text{Co}_3) = \text{Co}_3$ we end up with $A = G \cong \text{Co}_3$.
\end{proof}

The most delicate part in the previous proof is the computation of the right hand side of $\eqref{zaehlen}$. In the following we explain in greater detail this time consuming task (implementation in \texttt{PARI/GP} \cite{PARI2} with a total computing time of about 8 days using 550 threads simultaneously at the \textit{High Performance Computing Cluster} at the University of Würzburg).

For the sake of simplicity we write $f:=p_\mathfrak{p} - t_0q_\mathfrak{p}\in \mathbb{F}_\lambda [X]$ for some $t_0 \not\in S$.
In the case $k=3$ the following holds: If $f$ has exactly $d_i$ irreducible $\mathbb{F}_\lambda [X]$-factors of degree $i$ for $i\in \{1,2,3\}$ then $\pi_3(\text{Frob}(t_0)) = \binom{d_1}{3} + d_1d_2 + d_3$. Note that if a specialization reduces the degree, we have to add 1 to $d_1$.

In order to find $d_1$ we compute $p_1:= \gcd(X^\lambda-X,f)$. Clearly, $d_1 = \deg(p_1)$. Since $\lambda$ is too large for an efficient computation, we replace $X^\lambda-X$ with its reduction modulo $f$, which can be determined by the \textit{exponentiation by squaring}-method. 
In the same fashion we find $d_2$ and $d_3$: For $p_2 := \gcd(X^{\lambda^2}-X , \frac{f}{p_1})$ and $p_3:= \gcd(X^{\lambda^3}-X , \frac{f}{p_1p_2})$ we have $d_2 = \frac{1}{2}\deg(p_2)$ and $d_3 = \frac{1}{3}\deg(p_3)$.

Partial results for the computation of the right hand side of \eqref{zaehlen} can be found in the ancillary \texttt{Magma}-readable file.

\subsection{The symplectic group \boldm{$\text{PSp}_6(2)$}} \label{PSP}
In \cite[Theorem 4.2, Theorem 5.2]{BKW} both authors and Joachim König computed four-branch-point covers of degrees $28$ and $36$ with respective geometric monodromy group $G$ isomorphic to the $2$-transitive symplectic group $\text{PSp}_6(2)$. In order to verify $G=\text{PSp}_6(2)$, standard techniques yield that $G$ is either $\text{PSp}_6(2)$ or an alternating group. In contrast to the arguments given in \cite{BKW} to rule out the last case we now apply Elkies' method to give an alternative proof for $G = \text{PSp}_6(2)$. Assume, $G$ is $3$-transitive, then for the above covers we get a contradiction regarding the Hasse-Weil bound:
\begin{center}
\renewcommand{\arraystretch}{1.3}
\begin{tabular}{c||c|c}
degree &  28 & 36\\ \hline
ramification type & {\small$(2^6.1^{16},2^{12}.1^{4},2^{12}.1^{4},7^4)$} & {\small$(3^{12},2^{12}.1^{12}, 2^{12}.1^{12}, 4^7.2^1.1^6)$ }  \\\hline
$g(C_3)$ & 396 &  1275  \\\hline
$\lambda$ & 700\,001  & 7\,000\,003 \\\hline
Hasse-Weil bound & $\leq$ 1\,362\,637 & $\leq$ 13\,746\,671\\ \hline
$\#C_3(\mathbb{F}_\lambda)$ & $\geq$ 1\,405\,359  & $\geq$ 14\,032\,224   \\\hline
computing time & $\approx$ 2 minutes & $\approx$ 35 minutes
\end{tabular}
\end{center}

\subsection{The linear group \boldm{$\text{PSL}_6(2)$}} \label{PSL}
In \cite[Theorem 3.3]{BW} both authors calculated degree-$63$ four-branch-point covers with geometric monodromy group $G$ isomorphic to the $2$-transitive group $\text{PSL}_6(2)$.  Again, the main task in the verifying process boils down to deciding whether $G$ is either $\text{PSL}_6(2)$ or a full alternating/symmetric group. By applying Elkies' method we are able to give an alternative proof that $G$ cannot be $3$-transitive:
\begin{center}
\renewcommand{\arraystretch}{1.3}
\begin{tabular}{c||c}
degree &  63 \\ \hline
ramification type & {\small$(2^{28}.1^{7},  3^{20}.1^3, 3^{20}.1^3, 2^{16}.1^{31})$}    \\\hline
$g(C_3)$ & 5300  \\\hline
$\lambda$ & 120\,000\,007 \\\hline
Hasse-Weil bound & $\leq$ 236\,117\,193 \\ \hline
$\#C_3(\mathbb{F}_\lambda)$ & $\geq$ 239\,980\,524     \\\hline
computing time & $\approx$ 5 days
\end{tabular}
\end{center}

\subsection*{Computational remark}

In the accompanying file we provide a \texttt{Magma}-program to illustrate the computation of the right hand side of \eqref{zaehlen} for $M_{23}$, $\text{PSp}_6(2)$ and $\text{PSL}_6(2)$.

The specified computing times for these examples refer to computers with an \textit{AMD Ryzen 7 3700X} processor.

\section*{Acknowledgements}
We would like to thank Stephan Elsenhans and Peter Müller for some helpful discussions.

\end{document}